\documentclass[12pt,letterpaper,oneside]{article}
\usepackage{amsthm}
\usepackage[intlimits]{amsmath}
\usepackage{times}
\usepackage{graphics}
\usepackage{lipsum}
\usepackage{amsfonts}
\usepackage{amssymb}
\usepackage[intlimits]{amsmath}
\usepackage{epsfig}
\usepackage{microtype}
\newtheorem{thm}{Theorem}[section]
\newtheorem{Theorem}[thm]{Theorem}

\newtheorem{Remark}[thm]{Remark}

\newtheorem{Corollary}[thm]{Corollary}

\newcommand\ackname{Acknowledgments}
\if@titlepage
  \newenvironment{acknowledgements}{%
      \titlepage
      \null\vfil
      \@beginparpenalty\@lowpenalty
      \begin{center}%
        \bfseries \ackname
        \@endparpenalty\@M
      \end{center}}%
     {\par\vfil\null\endtitlepage}
\else
  
\fi
\makeatother

\title {A survey of the regular weighted Sturm-Liouville problem - The non-definite case \footnote{This research is supported, in part, by grant U0167 of the Natural Sciences and Engineering Research Council of Canada and a N.S.E.R.C. University Research Fellowship.} \\ PREPRINT\footnote{The original version of this article appeared with minor changes in ``A survey of the regular weighted Sturm-Liouville problem: the non-definite case", in {\it Applied Differential Equations}, World Scientific, Singapore and Philadelphia, 1986: 109-137. Some typographical errors have been corrected in this version.}}
\author{Angelo B. Mingarelli\footnote{Respectfully dedicated to my wife, Leslie Jean} \\ University of Ottawa \\ \\ Current address: School of Mathematics and Statistics, \\ Carleton University, Ottawa, ON.,\\ Canada K1S 5B6 \\  amingare@math.carleton.ca }
\begin{document}
\maketitle
\pagenumbering{roman}
\pagestyle{plain}
\pagenumbering{arabic}

\section{Introduction}
\renewcommand{\theequation}{\thesection.\arabic{equation}}

Let $p,q,r: \,[a,b]\;\rightarrow\;\mathbb{R}$ where $-\infty<a<b<\infty$ and $p(x)>0$ a.e., $q,r,\frac{1}{p}\,\in\,L(a,b)$ and
\begin{equation}
\label{eq1}
\int \limits_a^b |r(s)|\,ds\,>\,0.
\end{equation}
The weighted regular Sturm-Liouville problem consists in finding the values of a parameter $\lambda$ (generally complex) for which the equation
\begin{equation}
\label{eq2}
-(p(x)y^{\prime})^{\prime}+q(x)y=\lambda r(x)y \quad a\leq x\leq b,
\end{equation}
has a solution y (non-identically zero) satisfying a pair of homogeneous separated boundary conditions,
\begin{equation}
\label{eq3}
y(a)\cos \,\alpha\,-\,(py')(a)\sin\,\alpha=0
\end{equation}
\begin{equation}
\label{eq4}
y(b)\cos\,\beta\,+\,(py')(b)\sin\,\beta=0
\end{equation}
where $0\leq \alpha,\beta\,<\,\pi.$

Let $D\,=\,\{y:[a,b]\,\rightarrow\,\mathcal{C}\,|\, y, py'\,\in\,AC[a,b],r^{-1}\{(-py')'+qy\}\,\in\,L^2(a,b), y \;\text{satisfies} \eqref{eq3}-\eqref{eq4}\}$, in the case $|r(x)|>0$ a.e. on (a,b). Then associated with the problem \eqref{eq2}-\eqref{eq4} are the quadratic forms L and R, with domain D, where, for $y\,\in\,D$,
\begin{equation}
\label{eq5}
(Ly,y)=|y(a)|^2\cot\, \alpha\,+\,(py')(a)y'(a)\cot\,\beta\,+\,\int\limits_a^b\{p(x)|y'|^2\,+\,q(x)|y|^2\}\, dx
\end{equation}
and
\begin{equation}
\label{eq6}
(Ry,y)=\int\limits_a^b r(x)|y|^2\, dx.
\end{equation}
Here (,) denotes the usual $L^2$-inner product. (Moreover we note, as is usual, that the $\cot\,\alpha$ (resp. $\cot\,\beta$) term in \eqref{eq5} is absent if $\alpha=0$ (resp. $\beta=0$) in \eqref{eq3}-\eqref{eq4}).

It was noticed at the turn of this century by Otto Haupt \cite{twentynine},\cite{thirty} and Roland (George Dwight) Richardson, \cite{fiftysix}, \cite{fiftyseven}, \cite{sixtyone} (and it is very likely that David Hilbert was also aware, since Richardson went to G\"{o}ttingen in 1908...) that the nature of the boundary problem \eqref{eq2}-\eqref{eq4} is dependent upon some "definiteness" conditions on the forms L and R.

Thus, Hilbert and his school, termed the problem \eqref{eq2}-\eqref{eq4} \textbf{polar} if the form (Ly,y) is definite on D, i.e., either $(Ly,y)>0$ for each $y\neq 0$ in D or $(Ly,y)<0$ for each $y\neq 0$ in D (modern terminology refers to this case as the \textbf{left-definite} case, cf., \cite{twentytwo}).

The problem \eqref{eq2}-\eqref{eq4} was called \textbf{orthogonal} (or \textbf{right-definite} nowadays) if R is definite on D, (see above), whereas the general problem, i.e., when neither L nor R is definite on D, was dubbed \textbf{non-definite} by Richardson (see \cite{sixtyone}, p.285). In this respect, see also Haupt (\cite{thirty}, p.91).
We retain Richardson's terminology ``non-definite" in the sequel, in relation to the ``general" Sturm-Liouville boundary problem \eqref{eq2}-\eqref{eq4}. Thus, in the non-definite case, there exists functions y, z in D for which $(Ly,y)>0$ and $(Lz,z)<0$ and also, for a possibly different set of y, z, $(Ry,y)<0$ and $(Rz,z)>0$.

\section{The early theory of the non-definite case}
\setcounter{equation}{0}  

The first published results pertaining albeit indirectly, to the non-definite Sturm-Liouville problem \eqref{eq2}-\eqref{eq4} appear to be due to Emil Hilb \cite{thirtytwo}, (see also \cite{thirtythree}), who considered the single equation
\begin{equation}
\label{two1}
y''\,+\,(A\phi(x)\,+\,B)y\,=\,0,\quad 0\leq x\leq 1
\end{equation}
in the two parameters A, B and extended Klein's oscillation theorem to this case. In \eqref{two1} $\phi$ is real and continuous and may change sign in (0,1). Note that \eqref{two1} is a special case of \eqref{eq2}, with $\lambda\equiv A,\,p\equiv 1$, $r\equiv \phi,\;q\equiv -B$. In \cite{thirtytwo} Hilb noted the existence of parameter values (A,B) (corresponding to the Dirichlet problem) which gave rise to a ``minimum oscillation number" for corresponding eigenfunctions, this, being one of the characteristics of non-definite problems, (however this may also occur in orthogonal problems, see [\cite{sixtyone}, p.294]). Both Haupt (\cite{thirty}, p.69 and p.88) and Richardson (\cite{fiftynine}, p.289) cite Hilb's paper {\cite{thirtytwo} although it is to be noted that Maxime B\^{o}cher (\cite{eleven}, p.173) only cites Richardson, in this context, (and the citation is by way of an oral communication of Richardson to B\^{o}cher (\cite{eleven}, p.173, footnote) in his important survey of Sturm-Liouville theory up to about the year 1912. For a further up-date on the early developments of Sturmian theory see Richardson's review article (\cite{sixtytwo}, pp.110-111) of B\^{o}cher's classic book \cite{thirteen}, Lichtenstein's paper \cite{thirtynine} and the classic book by Ince (\cite{thirtyfour},\S\S 10-11). We note, however, that in a later paper \cite{twelve}, B\^{o}cher discovered Hilb's article \cite{thirtytwo}, (see \cite{twelve}, p.7, footnote).

Now it appears as if the first published theoretical investigations of the general non-definite problem (i.e., \eqref{eq2}) with q not identically a constant function on (a,b) are due to Haupt \cite{twentynine} in his dissertation, and Richardson \cite{fiftyeight}. In each one of these works one finds the basis for a beautiful extension of the Sturm Oscillation Theorem in relation to a non-definite problem in which the coefficients are assumed, \textbf{a priori}, continuous over the finite interval under consideration. A revised version of Haupt's dissertation (along with a corrected version of the said oscillation theorem) appeared much later, in 1915, (see \cite{thirty}). Richardson's version of the ``general Sturm oscillation theorem" first appeared in \cite{fiftynine} but the correct statement is to be found in (\cite{fiftynine}, errata) after a comment by George D. Birkhoff.

In Haupt's paper \cite{thirty} one finds many interesting results on the nature of non-definite problems, results with which are, unfortunately, not mentioned explicitly in Richardson \cite{sixtyone} a paper which deals essentially with ``new" oscillation theorems for \eqref{eq2} in the non-definite case (\cite{sixtyone}, $\S$4) and for the difficult case when q, in \eqref{eq2}, is allowed to vary \textbf{non-linearly} with $\lambda$, (see also McCrea-Newing \cite{fourtytwo}). However, Richardson does refer to Haupt \cite{thirty} insofar as the oscillation theorem is concerned, as Haupt's theorem appeared one year earlier (cf., above). Thus the first version of an oscillation theorem (in the case q$\neq$ constant) is due to Haupt \cite{thirty}. On the other hand, a non-trivial sharpening, and the final beautiful form of the oscillation theorem, was formulated by Richardson in (\cite{fiftynine}, errata), (see also Richardson (\cite{sixtyone}, p.285)). Curiously enough Richardson was indeed the ``father of the non-definite case" as one can infer from his remarks in (\cite{sixtyone}, p.285), the claim being that the first investigations of a non-definite problem appeared in Richardson \cite{fiftyeight}.

It is remarkable that Haupt's work \cite{thirty} has remained obscure, even today. For example Ince (\cite{thirtyfour}, p.248) refers to Richardson but nowhere in the book does one find a reference to Haupt! One reason for this may be the following: As was mentioned above Hilbert and his school referred to left-definite problems as ``polar" problems and to the right-definite as ``orthogonal" problems. However, Haupt, in \cite{thirty}, decided to include ``non-definite" problems (in Richardson's terminology) under the heading of ``polar" ones, Haupt (\cite{thirty}; p.84, $\S$5), whereas the ``orthogonal" case was included in the section entitled ``non-polar case", Haupt (\cite{thirty}, p.76, $\S$3). It may then have appeared, to those influenced by Hilbert and his school, that, at first sight, Haupt was not doing anything radically new except for extensions, to equations in which the parameter $\lambda$ appears non-linearly, of the Sturm oscillation theorem, Haupt (\cite{thirty}, $\S$4). Haupt's failure to emphasize the \textbf{"non-definite polar case"}, Haupt (\cite{thirty}, p.91), as a new and distinct case may have led to a low readership of his paper, \cite{thirty}.

\begin{Theorem}\emph{Haupt (\cite{thirty}, pp.84-85)}\\
\label{thm1}
In \eqref{eq2}, let $p, q, r$ be continuous in $[a,b]$ along with $p$ differentiable and $p(x)>0$ there. Then the eigenvalues of \eqref{eq2}-\eqref{eq4} are the zeros of an entire transcendental function which is not identically zero.
\end{Theorem}
\begin{Remark}
For results pertaining to the \textbf{order} of the entire transcendental function alluded to by Haupt, see Halvorsen \cite{twentyseven}, Mingarelli (\cite{fourtyeight,fifty,fourtyseven}) and Atkinson-Mingarelli \cite{four}.
\end{Remark}

Of course, it is necessary to assume in Theorem~\ref{thm1} that there hold \eqref{eq1} or else the entire function mentioned therein may, in fact, vanish identically on $\mathcal{C}$ as an easy example will show. This hypothesis is not explicit in Haupt \cite{thirty} but is, however, necessary in the proof.

\begin{Corollary}\label{cor1}
The eigenvalues of \eqref{eq2}-\eqref{eq4} form a discrete subset of the complex plane (i.e., having no finite point of accumulation.)
\end{Corollary}

In Haupt \cite{thirty} a ``normalization" is assumed for \eqref{eq2}, i.e., the equation is transformed into a ``normal form" in many places (see Haupt (\cite{thirty}, pp.85-86) for the transformation and Mingarelli \cite{fourtyfour} for an assumption in the same vein).

The existence of possibly non-real eigenvalues was alluded to by Haupt (\cite{thirty}, pp.94) and by Richardson (\cite{sixtyone}, p.289 and footnote) however neither author gave an instance of such an occurrence (cf., also Mingarelli (\cite{fourtysix}, Chapter 4)). Concrete examples of non-definite problems with non-real eigenvalues were obtained in Mingarelli (\cite{fourtyfour}, pp.519-520), (\cite{fourtyseven}, p.376) in which some very early ideas of Hilb \cite{thirtythree} were used to show that $\lambda=i$ may be turned into an eigenvalue for a Dirichlet problem associated with \eqref{eq2}, (see \cite{fourtyseven} for details).

\section{Terminology and notation}
The theory of linear operators in spaces with an indefinite metric, Azizov and Iohvidov \cite{five}, Bogn\'{a}r \cite{fourteen}, together with their applications to quantum field theory, Nagy \cite{fiftythree}, suggests the following terminology (much of which was due to Werner Heisenberg).

An eigenfunction $y$ of \eqref{eq2}-\eqref{eq4} corresponding to a non-real eigenvalue will be called a \textbf{complex ghost state} or complex ghost. (These necessarily satisfy $\int\limits_a^b r(x)|y|^2\, dx=0$.)

If the non-real eigenvalue $\lambda$ is non-simple its corresponding eigenfunction y will be said to be \textbf{degenerate}. It is \textbf{non-degenerate} otherwise, i.e., $\int\limits_a^b r(x)|y|^2\, dx\neq0$, (see $\S$5).

An eigenfunction corresponding to a real non-simple eigenvalue will be called a \textbf{degenerate real ghost state} (or a \textbf{dipole ghost}, see \cite{fiftythree}) whereas a real eigenfunction y corresponding to the (real) eigenvalue $\lambda$ is a \textbf{non-degenerate real ghost state} provided $$\text{sign}\bigg\{\lambda \int\limits_a^b r(x)|y|^2\, dx\bigg\}<0.$$

The term \textbf{positive eigenfunction} (or \textbf{ground state}) will refer to a (real) eigenfunction whose values are strictly positive in the interior of the interval under consideration.

If $y$ is any eigenfunction the quantity $$\int\limits_a^b r(x)|y|^2\, dx$$ will be called the \textbf{$r$-Kre\v{i}n norm} of $y$. (This terminology is motivated by the formalism whereby one may cast the setting for a spectral theory of non-definite metric, e.g., Kre\v{i}n space, Pontryagin space. (see, e.g., Azizov and Iohvidov \cite{fourteen}, Bogn\v{a}r \cite{five} and applications in Langer \cite{thirtyseven}, Daho and Langer \cite{sixteen, seventeen}, \v{C}urgus and Langer \cite{fifteen}, Daho \cite{eighteen, nineteen}, and Mingarelli \cite{fourtysix}). Note that when $r(x)>0$ a.e. on $(a,b)$ the $r$-Krein-norm of $y$ is the usual norm in the weighted Hilbert space $L_r^2(a,b)$.

As is customary we will denote by $r_{\pm}$ the positive (negative) part of r, i.e., $r_{\pm}\,= \,max\{\pm r(x),0\}.$

\section{Some motivation}
Consider the boundary problem associated with Mathieu's equation
$$y''+(-\alpha+\beta p(x))y=0,\;-\infty<x<+\infty$$ where $(\alpha,\beta)\,\in\,{\Re}^2$ are parameters and $p$ is periodic, or more generally, locally Lebesgue integrable over $\Re$. In the usual problem, $\beta$ is fixed and values of $\alpha$ are sought $\ldots$. However if we fix $\alpha<<0$ and seek $\beta$'s the problem is generally non-definite. (We recall that Hilb \cite{thirtytwo} considered the same problem, but on a finite interval.)

More generally, the equation, $$y''+(-\alpha A(x)+\beta B(x))y=0,\;-\infty<x<+\infty,$$ in the two parameters $\alpha,\beta\,\in\,L_{loc}$($\Re$ ) also leads to a non-definite problem once one of the parameters $\alpha,\beta$ is held fixed and the other is treated as an eigenvalue parameter. The question of the existence of positive solutions of the latter equation has received some interest lately (e.g., the monograph of Halvorsen-Mingarelli \cite{twentyeight}). For applications of these to laser theory see Heading \cite{thirtyone}, also McCrea-Newing \cite{fourtytwo} and the references therein.

Techniques from the theory of non-definite problems were recently used by Deift-Hempel \cite{twenty} with applications to the theory of color in crystals. Various other related problems are treated, for instance, in Barkovskii-Yudovich \cite{six,seven} in relation to Taylor vortex formation arising from rotating cylinders.

\section{On the existence of eigenvalues}

The problem of the actual existence of eigenvalues for the non-definite problem \eqref{eq2}-\eqref{eq4} is treated implicitly in Haupt \cite{thirty} and somewhat more elaborately in Richardson \cite{sixtyone}, at least for the Dirichlet problem. Pr\"{u}fer angle methods have been used in this respect in Halvorsen \cite{twentyseven} to, at least, settle the existence of an infinite sequence of \textbf{real} eigenvalues under the more general set of boundary conditions \eqref{eq3}-\eqref{eq4}. A modified Pr\"{u}fer angle method has been used in Atkinson-Mingarelli \cite{four} which also settles the existence and, at the same time, yields their asymptotics.
\begin{Theorem}\emph {Haupt \cite{thirty}, Richardson \cite{sixtyone}, Halvorsen \cite{twentyseven}, Atkinson-Mingarelli \cite{four}}.\\
\label{thm51}
\noindent Whenever $$\int\limits_a^b r_{+}(s)ds>0\;\; \text{and}\;\;\int\limits_a^b r_{-}(s)\,ds>0$$ the problem \eqref{eq2}-\eqref{eq4} has two infinite sequences of real eigenvalues, one positive and one negative, and each one of which has $+\infty$ and $-\infty$ for its only point of accumulation.
\end{Theorem}
\begin{Remark}
In Haupt \cite{thirty}, Richardson \cite{sixtyone} the authors deal with the case in which the coefficients $p,q,r$ are continuous in $[a,b]$ whereas in Halvorsen \cite{twentyseven} and Atkinson-Mingarelli \cite{four} the more general case referred to in the introduction is considered.
\end{Remark}
For a non-definite problem \eqref{eq2}-\eqref{eq4}, non-real eigenvalues may or may not occur! Both these cases are possible - with regards to the former see the comments in Haupt \cite{thirty}, Richardson \cite{sixtyone}, Daho and Langer \cite{sixteen, seventeen} and the specific example in Mingarelli (\cite{fourtyfour}, p.250) or Atkinson-Jabon (\cite{three}, appendix), Jabon \cite{thirtyfive}, or Marziali \cite{fourtyone}.

The question of the possible existence of non-real eigenvalues for a non-definite problem \eqref{eq2}-\eqref{eq4} seems to have first been formulated by Haupt \cite{thirty} and then by Richardson (\cite{sixtyone}, p.100) for the case of continuous coefficients under the assumption that the parameter is ``normalized" (see Haupt \cite{thirty}, pp.85-86). A result in the same vein was obtained independently by the author in Mingarelli (\cite{fourtyfour}, Theorem 2) the proof of which also includes the case of Lebesgue integrable coefficients.

\begin{Theorem}\emph{(Haupt \cite{thirty}, Mingarelli \cite{fourtyfour})}\\
\label{thm52}
\noindent In the non-definite case and under the assumptions of $\S$, the problem \eqref{eq2}-\eqref{eq4} has at most finitely many non-real eigenvalues, their total number being even.
\end{Theorem}

A particular case of a non-definite case of a non-definite (singular) problem was explored in Daho-Langer \cite{sixteen, seventeen}, (see also \v{C}urgus-Langer \cite{fifteen} and Mingarelli \cite{fourtyfour}, Chapters 3,4).

There is indeed an intimate connection between the theory of symmetric linear operators in a Pontryagin space (Krein space) and the theory of non-definite Sturm-Liouville problems. We will not delve into this matter here, for the sake of brevity, although we will refer the reader to the relevant literature herewith: For general notions on Pontryagin spaces see the monograph by Bogn\'{a}r \cite{fourteen} and the excellent survey paper by Azizov-Iohvidov \cite{five}. For applications of this theory to the problem at hand see Langer \cite{thirtyseven}, Daho-Langer \cite{sixteen, seventeen}, \v{C}urgus-Langer \cite{fifteen}, while for interesting physico-theoretical applications of Pontryagin spaces we refer to the monograph by Nagy \cite{fiftythree}. We need only mention at this point that Theorem~\ref{thm52} and other similar theorems, some of which we will present below, may be proved with the help of the theory of symmetric linear operators in a Pontryagin space.

Now for $\lambda$ real let n($\lambda$) denote the number of negative eigenvalues (counting multiplicities) of the problem
$$-(p(x)y')'+(q(x)-\lambda r(x))y=\nu y$$ where $y$ is required to satisfy \eqref{eq3}-\eqref{eq4}. Then the function n($\lambda$) has an absolute minimum which will be denoted by $n_0$. (This construction is due to Haupt (\cite{thirty}, p.85.)
\begin{Theorem}\emph{Haupt \cite{thirty}, Mingarelli \cite{fourtyfour}}\\
\label{thm53}
\begin{enumerate}
\item Let $p,q,r$ be continuous on $[a,b]$ and \eqref{eq2}-\eqref{eq4} non-definite. Then the number of pairs (i.e., an eigenvalue and its complex conjugate) of non-real eigenvalues does not exceed $n_0$, whenever the parameter is normalized, (\cite{thirty}, p.100).
\item In (i) above we may replace $n_0$ by the number of negative eigenvalues of the problem
\begin{equation}\label{eq51}
-(p(x)y')'+(q(x)-\lambda r(x))y=\lambda y
\end{equation}
subject to \eqref{eq3}-\eqref{eq4} provided zero is not an eigenvalue of the said problem, \cite{fourtyfour}. Cases of equality here may be exhibited, see \cite{fourtyfour} or Marziali \cite{fourtyone}.
\end{enumerate}
\end{Theorem}

We now note that, using a Green's function argument, in the regular case of the non-definite problem \eqref{eq2}-\eqref{eq4} the spectrum is purely discrete (i.e., it consists only of eigenvalues of finite multiplicity, Mingarelli \cite{fourtythree}.) We now turn our attention to the question of the ``simplicity" of the eigenvalues. For basic results regarding the nature of ``simple" and ``non-simple" eigenvalues we refer to Ince (\cite{thirtyfour}, $\S$ 10.72). Thus let $y(x,\lambda)$ be a non-trivial solution of \eqref{eq2} which satisfies \eqref{eq3}. Then the function
$$F(\lambda)\equiv y(b,\lambda)\cos\,\beta+(py')(b,\lambda)\sin\,\beta$$ is an entire function of $\lambda \in\mathcal{C}$ whose order, generally, does not exceed one-half (see e.g., Atkinson \cite{one}). Actually it is now known that its order is \textbf{precisely} one-half as was shown directly by Halvorsen \cite{twentyseven} and, as a consequence of the asymptotics, by Mingarelli \cite{fourtyseven}, (cf., Mingarelli \cite{fourtyeight}), and Atkinson-Mingarelli \cite{four}. For some extensions of the above result on the order of $F$ see Mingarelli \cite{fifty}.

Now the zeros of $F$ are in a one-to-one correspondence with the eigenvalues of the problem \eqref{eq2}-\eqref{eq4}. We say that an eigenvalue $\lambda\,\in\,\mathcal{C}$ is \textbf{simple} if it is a simple zero of $F$ (i.e., $F'(\lambda)\neq0$). It is said to be \textbf{non-simple} otherwise, (i.e., if $0=F(\lambda)=F'(\lambda)$).

\textbf{Note}: It does not follow from the above definition that a non-simple eigenvalue necessarily has two linearly independent eigenfunctions associated with it, for if this were ever the case, their span would generate the space of \textbf{all} solutions of \eqref{eq2} and clearly one could find solutions which do not satisfy the first boundary condition \eqref{eq3}, (cf., Ince \cite{thirtyfour}, p.241).

Thus \textbf{for the problem under consideration, namely \eqref{eq2}-\eqref{eq4}, there is a one-to-one correspondence between any real eigenvalue and its corresponding eigenfunction} (if it is suitably normalized).

It is known that $\lambda\,\in\,\mathcal{C}$ is \textbf{non-simple if and only if it has a corresponding (real or complex) eigenfunction which is a degenerate ghost state}. The latter result was anticipated by Richardson (\cite{sixtyone}, p.294, Corollary) albeit without proof. For a proof and an extension of the said result to the case of measurable coefficients, see Mingarelli \cite{fourtynine}.

It turns out that \textbf{real degenerate ghost states may exist for \eqref{eq2}-\eqref{eq4}}, see the example in Mingarelli \cite{fourtynine}. Furthermore there \textbf{may also exist non-degenerate real ghost states}, see, once again, an example in Mingarelli \cite{fourtynine}, thus answering a question raised by Haupt (\cite{thirty}, p.100 footnote). For all of the examples now known it appears that the complex ghost states are all non-degenerate. \textbf{There is no known example of a degenerate complex ghost state} although the author feels that these very likely do exist in some cases.
\begin{Theorem}\emph{Mingarelli \cite{fourtynine}}\\
\label{thm54}
\noindent Let $\lambda>0$ and assume that zero is not an eigenvalue of \eqref{eq2}-\eqref{eq4}. Then the total number of non-degenerate and degenerate real ghost states is always finite and bounded above by the number of negative eigenvalues of \eqref{eq51}-\eqref{eq3}-\eqref{eq4}. An analogous result holds for $\lambda<0$.
\end{Theorem}

\begin{Remark}
Specific examples, Marziali \cite{fourtyone}, seem to indicate that the bound appearing in Theorem \ref{thm54} is sharp. Consolidating the above results we may formulate,
\end{Remark}

\begin{Theorem}
In the non-definite case of \eqref{eq2}-\eqref{eq4} the spectrum is discrete, always consists of a doubly infinite sequence of real eigenvalues, having no finite limit, and has at most a finite and even number of non-real eigenvalues (necessarily occurring in complex conjugate pairs) along with at most finitely many real non-simple eigenvalues. The totality of all such eigenvalues comprise the spectrum of \eqref{eq2}-\eqref{eq4}.
\end{Theorem}

For an extension of Theorem~\ref{thm53}-\ref{thm54} to an abstract setting, see Mingarelli \cite{fiftyone} wherein the more general operator equation $Ax=\lambda Bx$ is considered as a generalized eigenvalue problem in a complex Hilbert space, thus allowing for extensions of the said results to partial differential equations with indefinite weight-functions.

\textit{Note}  \textbf{Open Problems}
\begin{enumerate}
\item For a given non-definite problem \eqref{eq2}-\eqref{eq4} can one find an \textbf{a priori} bound on the modulus (or real/imaginary part) of the ``largest" non-real eigenvalue which might appear?
\item For a given non-definite problem \eqref{eq2}-\eqref{eq4} can one find an \textbf{a priori} bound on that interval of the real axis which contains all the non-simple eigenvalues and those eigenvalues corresponding to non-degenerate real ghost states? (Note that this interval is finite on account of Theorem~\ref{thm54}). In relation to this question see Atkinson-Jabon \cite{three} wherein this is done for a specific example.
\item Find a \textbf{sufficient} condition which will guarantee the existence of a non-real eigenvalue for a non-definite problem \eqref{eq2}-\eqref{eq4}. (\textbf{Necessary} conditions are widespread: see e.g. Mingarelli (\cite{fourtyfour}, p.521, Theorem 1.)
\item According to Richardson (\cite{sixtyone}, p.289) all sufficiently large real eigenvalues may be regarded as furnishing a minimum of a calculus of variations problem (since they all have positive r-Krein-norm?); however, no proof of this claim is given.
\end{enumerate}

(Note, however, one of Richardson's last papers on this subject, Richardson \cite{sixtythree}.)

\section{A generalized Sturm Oscillation Theorem}
The classical Sturm Oscillation Theorem for a left-definite Dirichlet problem associated with \eqref{eq2} states that, if $r$ satisfies the hypothesis of Theorem~\ref{thm51}, the n-th positive (negative) eigenvalue has an eigenfunction which has precisely $n$ zeros in $(a,b)$, (see e.g. Ince \cite{thirtyfour}, p.235, Theorem 3).

\begin{Theorem}\emph{Haupt \cite{thirty}, Richardson \cite{sixtyone}}
\label{thm61}
In the non-definite case of \eqref{eq2}-\eqref{eq4} there exists an integer $n_R\geq0$ such that for each $n\geq n_R$ there are \textbf{at least} two solutions of \eqref{eq2}-\eqref{eq4} having exactly $n$ zeros in $(a,b)$ while for $n< n_R$ there are \textbf{no real solutions} having $n$ zeros in $(a,b)$. Furthermore there exists a possibly different integer $n_H\geq n_R$ such that for each $n\geq n_H$  there are precisely two solutions of \eqref{eq2}-\eqref{eq4} having exactly $n$ zeros in $(a,b)$.
\end{Theorem}
\begin{Remark}
We will call $n_R$ the \textbf{Richardson Index} {or number} and label $n_H$ the \textbf{Haupt Index} (or number) of the problem \eqref{eq2}-\eqref{eq4} for historical reasons- The existence of $n_H$ appears to have been first established in the literature by Haupt (\cite{thirty}, p.86) while the existence of $n_R$, in general, is almost certainly due to Richardson \cite{fiftyeight} and \cite{sixtyone}, (modulo Hilb's special case \cite{thirtytwo}).
\end{Remark}

It now follows from the \textbf{Haupt-Richardson oscillation theorem} (Theorem~\ref{thm61}) that, generally speaking, a non-definite problem will tend not to have a real ground state (positive eigenfunction)!

The existence of $n_R$ in the case of measurable coefficients was obtained by the author, Mingarelli \cite{fiftytwo}, via a Pr\"{u}fer transformation and use of the results in Mingarelli \cite{fourtynine}.

It was shown in Mingarelli \cite{fourtyfive} that \textbf{each one of the cases $n_H=n_R$ and $n_H>n_R$ may occur}: Indeed the problem \eqref{eq2} with $p\equiv 1,\;q\equiv -9\pi^2/16,\;r(x)=+1$ on $[0,1]$ and $r(x)=-1$ on $(1,2]$, $\alpha=\beta=0$ on $[a,b]=[0,2]$ has precisely one pair of non-real eigenvalues (actually pure imaginary) situated at around $\pm 4.3$i while the remaining real eigenvalues have eigenfunctions with at least \textbf{one} zero in $(a,b)$, so that $n_R\geq 1$ and in fact $n_R=n_H=1$, in this case.

On the other hand if we set $q\equiv -22.206$ (a lower approximation to $-9\pi^2/4$) and define the other quantities as in the above example, an explicit calculation shows that $n_H=3$ while $n_R=2$. More recently the calculations reported in Atkinson-Jabon \cite{three} and Marziali \cite{fourtyone} yield essentially unlimited examples of such pathological behavior (see also Richardson \cite{sixtyone}, p.298).

Armed with the Haupt-Richardson oscillation theorem one may proceed to prove asymptotic estimates for the real eigenvalues of non-definite problems as was done in Mingarelli \cite{fourtyfour}- This is closely related to a conjecture of the late Konrad J\"{o}rgens \cite{thirtysix}. Under assumptions which, in the aftermath, are very similar to those used by M.H. Stone \cite{sixtyfour}, (e.g., if $r$ is continuous, then $r(x)$ changes sign finitely many times in $(a,b)$), the author showed

\begin{Theorem} (\emph{Mingarelli \cite{fourtyseven}}.)
\label{thm62}
If the positive eigenvalues $\lambda_n^+$ of a given non-definite problem \eqref{eq2}-\eqref{eq4} are labeled according to the Haupt-Richardson oscillation theorem so that $\lambda_n^+$ has an eigenfunction with precisely $n$ zeros in $(a,b)$, then
\begin{equation}
\label{eq61}
\lambda_n^+ \sim \frac{n^2\pi^2}{\{\int_a^b \sqrt{(\frac{{r}}{p})_+}\, ds\}^2}, \quad n\to \infty,
\end{equation}
where $(\frac{{r}}{p})_+ = max\{(\frac{{r}}{p})_+,0\}$. An analogous formula holds for the negative eigenvalues as well.
\end{Theorem}

Let \textbf{N($\lambda$) denote the number of zeros of a non-trivial solution of \eqref{eq2} satisfying, say, \eqref{eq3}}. If we label our positive eigenvalues so that $N(\lambda_n^+)=n$ for $n\geq n_H$, then it will follow from \eqref{eq61} that
\begin{equation}
\label{eq62}
N(\lambda)\sim \sqrt{\lambda}\pi^{-1}\int_a^b \sqrt{(\frac{{r}}{p})_+}\,ds,\quad \lambda\to \infty,
\end{equation}
This last relationship entails what the author calls \textbf{J\"{o}rgens' conjecture} (\cite{thirtysix}, p.5.16). A proof of this conjecture, in general, is given, among other results, in Atkinson-Mingarelli \cite{four} thus settling the question of the validity of \eqref{eq62}, at least in the case when $p(x)>0$ a.e. on $(a,b)$. Now let \textbf{$n_+(\lambda$) denote the number of positive eigenvalues of \eqref{eq2}-\eqref{eq4} which are in $(0,\lambda)$}. Then it follows from  \eqref{eq61} once again that
\begin{equation}
\label{eq63}
n_+(\lambda)\sim \sqrt{\lambda}\pi^{-1}\int_a^b \sqrt{(\frac{{r}}{p})_+}\,ds,\quad \lambda\to \infty,
\end{equation}
(for a detailed calculation see Mampitiya (\cite{fourty}, p.28, Lemma 4.4). We mention, in passing, that \eqref{eq63} has been recently extended to incorporate polar \textbf{vector} Sturm-Liouville problems on a finite interval, see Mampitiya \cite{fourty}.

One of the truly remarkable results that one may extract from Richardson \cite{sixtyone} is the following.

\begin{Theorem} (\emph{\cite{sixtyone}, p.302, Theorem 10}.)
\label{thm63}
Let $r$ be continuous and not vanish identically in any right-neighborhood of $x=a$. If $r(x)$ changes its sign precisely once in $(a,b)$ then the roots of the real and imaginary parts $u,v$, of any non-real eigenfunction $y=u+iv$ corresponding to a non-real eigenvalue, separate one another.
\end{Theorem}
\begin{Corollary}\label{cor64}
Under the conditions of Theorem~\ref{thm63} it follows that any non-real eigenfunction $y$ of \eqref{eq2}-\eqref{eq4} cannot have a zero for $x$ in $(a,b)$.
\end{Corollary}

\textit{Note}   The initial assumption on $r$ in Theorem~\ref{thm63} is not explicitly mentioned in the proof given by Richardson \cite{sixtyone} however it appears to be necessary for the validity of his proof.

For an extension of Theorem~\ref{thm63} to the case when $r(x)$ changes sign its sign finitely many times, see Mingarelli (\cite{fourtyfour}, p.525, Theorem 3).

Finally we note that some interesting comparison theorems for non-definite problems are derived in Atkinson-Jabon (\cite{three}, p.35, Proposition 3.4).

\section{Further Results and Extensions}

With regards to an eigenfunction expansion in the regular non-definite case see Faierman \cite{twentythree} where the actual uniform convergence of such an expansion is also treated. Another approach to the same problem is to be found in Binding-Seddighi \cite{nine} with an aim towards an abstract expansion theorem which includes the non-definite case as a particular case, but with additional boundedness assumptions on the coefficients.

The problem of finding asymptotic expressions for the solutions of equations of the form \eqref{eq2} with additional smoothness on the coefficients, is an old one and much has been done in this direction- The book by Olver \cite{fiftyfour} is a good introduction to the subject as is the book by Erdelyi \cite{twentyone}, especially Chapter 4. Further results may be found in recent papers by Fedoryuk \cite{twentyfour,twentyfive} and the paper by Olver \cite{fiftyfive}, see also the references therein. Interesting results in this connection are also buried in Birkhoff-Langer \cite{ten}.

Results in the case of \textbf{singular} non-definite Sturm-Liouville problem (dealing almost exclusively with a closed half-line) are more or less scattered, although the object of recent interest. In this direction see the pioneering papers by Atkinson-Everitt-Ong \cite{two}, and Daho-Langer \cite{sixteen,seventeen} wherein an expansion theorem is formulated. The paper by Langer \cite{thirtyseven} sheds much insight into the role that Krein spaces play in the study of non-definite problems. Further expansion results (full and half-range expansions) for higher-order scalar ordinary differential equations are given in \v{C}urgus-Langer \cite{fifteen}, see also Daho \cite{eighteen,nineteen} for results dealing with the existence of a Titchmarsh-Weyl matrix function in a higher-order singular case.

For the most up-to-date results regarding left-and right-definite problems, which include \eqref{eq2}, we refer to the paper by Everitt \cite{twentytwo} and Bennewitz-Everitt \cite{eight} and the references therein.

The basis for the extension of the foregoing results to partial differential equations may be found in Fleckinger-Mingarelli \cite{twentysix}, see also Mingarelli \cite{fiftyone}.

At this point many questions even in the regular case, remain unanswered: For example, questions relating to the existence/non-existence of real ground states for \eqref{eq2}-\eqref{eq4}; extensions of the preceding results to the systems of second order differential equations and higher-dimensional analogs; an efficient method for calculating the non-real eigenvalues of non-definite problems \eqref{eq2}-\eqref{eq4}, (current techniques of, Marziali \cite{fourtyone}, Atkinson-Jabon \cite{three}, rely upon explicit computation of the zeros of F($\lambda$), (see \S 5), which is generally inefficient; and to find \textbf{a priori} estimates on the Haupt and Richardson indices for a given non-definite problem.

\begin{center}

{\Large \bf Acknowledgments}
\end{center}
I am grateful to Professor F.V. Atkinson who, in 1977 introduced the author to this fascinating field. The main reference to Richardson \cite{sixtyone} was only discovered in 1978 whereas the importance of Haupt's work was only recognized in 1984. All the other references to the early history of the non-definite problems are basically recent finds. It is uncanny that so much work had been done in connection with this problem at the turn of this century, work which seemingly came to a halt around 1920. Renewed theoretical interest in this area appears to be very recent, only after a gap of almost 60 years!

I am also grateful to Dr. Philip Hartman (formerly of Johns Hopkins University) who, in 1978, remarked to me that Hilb may have had something to do with this problem. This hunch proved correct as we have seen. I am grateful to Phil Hartman for also raising the question, in 1980, of the ``simplicity" of the real eigenvalues of non-definite problems, a question which led to my paper \cite{fourtynine}.

I wish to thank the following mathematicians for supplying me with preprints of their work - C. Bennewitz, Paul Binding, Karim Daho, Percy Deift, Mel Faierman, and David Jabon. I also wish to thank Derick Atkinson, Paul Binding, Patrick Browne, Mel Faierman and Gotskalk Halvorsen for interesting discussions. My gratitude also goes to Hans Kaper, for providing the reference \cite{twentyone}, and Hubert Kalf for rekindling my interest in Lichtenstein's paper \cite{thirtynine}.

{\bf Note added June 29, 2011} The papers referred to as \cite{four}, \cite{nine}, \cite{twenty} and \cite{twentythree} below have since appeared and the updated citations are included in the references below. In addition, the papers by Langer \cite{thirtyeight}, Richardson \cite{sixty}, and  Wheeler \cite{sixtyfive} were inadvertently left out of the body of the original paper as it appeared 25 years ago. In \cite{thirtyeight} Langer makes a detailed study of the spectral theory of a regular Sturm-Liouville problem having only one turning point (of any order) on a finite interval. Wheeler \cite{sixtyfive} considered the existence and completeness questions in both the orthogonal and polar cases, while Richardson \cite{sixty} considers the polar and orthogonal case of a two dimensional elliptic problem while hinting at non definite cases as well (this seems to be the first case of such a study for higher dimensions), cf., [\cite{sixty}, p.494].

\bibliographystyle{plain}
\bibliography{Mingrefs2}

Department of Mathematics\\
University of Ottawa\\
Ottawa, Ontario, Canada\\
K1N 6N5
\end{document}